\documentclass[12pt,english]{article}
\usepackage{amsfonts,here,graphicx}
\usepackage{cite,amsfonts,amssymb} 
\def\qed{\kern8pt \vrule height5pt depth0pt width5pt}

\def\sup{\mathop{\rm sup}}

\def\proof{\noindent \medskip {\bf Proof:}$\;\;$}

\def\ass#1#2\endass{\vskip5pt plus2pt \noindent{\bf (A.#1)} {\it #2} \vskip5pt
plus2pt }

\newtheorem{prop}{Proposition}

\newtheorem{lem}{Lemma}

\newtheorem{theor}{Theorem}

\newtheorem{cor}{Corollary}

\newtheorem{defin}{Definition}

\title{\bf  Properties of convolutions arising in 
stochastic Volterra equations}

\author{\large\sf Anna Karczewska \\
 \\
 Department of Mathematics,
 University of Zielona G\'ora\\
 ul. Szafrana 4a, 65-246 Zielona G\'ora, Poland\\
 e-mail: A.Karczewska@im.uz.zgora.pl\\
}
\date{\today}

\begin{document}

\maketitle

 \noindent\def\thefootnote{}

\begin{abstract}
The aim of this note is to provide some results for
stochastic convolutions corresponding to stochastic Volterra equations
in separable Hilbert space. We study convolution of the form 
$W^{\Psi}(t):=\int_0^t S(t-\tau)\Psi(\tau)dW(\tau)$,
$t\geq 0$, where $S(t),~t\geq 0$, is so-called {\em resolvent} for 
Volterra equation considered,$\Psi$ is an appropriate process and $W$ is a 
cylindrical Wiener process. 
\end{abstract}

\vspace{8mm}
\begin{center}
{\large \em Abbr.\ title: \bf Stochastic Volterra Convolutions}
\end{center}

\newpage

\section{Definitions and notation} \label{defnot}

In the paper we consider the following stochastic Volterra equation 
in a separable Hilbert space $H$:
\begin{equation} \label{deq1}
 X(t) = X_0 + \int_0^t a(t-\tau) AX(\tau)\, d\tau 
 + \int_0^t \Psi(\tau)\,dW(\tau),
\end{equation}
where $t\in\mathbb{R}_+$, 
$a\in L_{\mathrm{loc}}^1(\mathbb{R}_+)$, $A$ is a closed 
unbounded linear operator in $H$ with a dense domain D(A),
$\Psi$ is an adapted integrable stochastic process specified below, $W$ is a
cylindrical Wiener process with respect to $t$ and $X_0$ belongs to $H$.

The stochastic Volterra equations have been studied in connection with 
applications to problems arising in mathematical physics, particularly
in viscoelasticity and heat conduction in materials with memory.
We refer to the papers  \cite{ClDaPr1},\cite{ClDaPr2} and \cite{ClDaPP}.
Let us note that the equation (\ref{deq1}) is a generalization of 
stochastic heat and wave equations and stochastic linear Navier-Stokes
system.

The above equation (\ref{deq1}) is a stochastic version of the deterministic 
Volterra equation of the form
\begin{equation} \label{deq2}
 X(t) = X_0 + \int_0^t a(t-\tau) AX(\tau) d\tau 
 + f(t),
\end{equation}
where elements in (\ref{deq2}) are the same as in (\ref{deq1}), and $f$ is an
appropriate $H$-valued mapping.

By $S(t),~t\geq 0$, we shall denote the family of resolvent 
operators corresponding to
the Volterra equation (\ref{deq2}) and defined as follows.

\begin{defin} \label{def1} (see, e.g.\ \cite{Pr})\\
A family $(S(t))_{t\geq 0}$ of bounded linear operators in the space $H$ is
called {\tt resolvent} for (\ref{deq2}) if the following conditions are satisfied:
\begin{enumerate}
\item $S(t)$ is strongly continuous on $\mathbb{R}_+$ and $S(0)=I$;
\item $S(t)$ commutes with the operator $A$, that is, $S(t)(D(A))\subset D(A)$
 and $AS(t)x=S(t)Ax$ for all $x\in D(A)$ and $t\geq 0$;
 \item the following {\tt resolvent equation} holds
\begin{equation} \label{deq3}
 S(t)x = x + \int_0^t a(t-\tau) AS(\tau)x d\tau 
\end{equation}
for all $x\in D(A),~t\geq 0$.
\end{enumerate}
\end{defin}
\noindent{\bf Comment:} Let us emphasize that the resolvent $S(t),~t\geq 0$,
is determined by the operator $A$ and the function $a$.
Moreover, as a consequence of the strong continuity of $S(t)$ we have
for any $T>0$
\begin{equation} \label{deq3a}
 \sup_{t\leq T} ||S(t)||<+\infty \;,
\end{equation}
where $||\cdot ||$ denotes the operator norm.

In the paper we shall assume that the equation (\ref{deq2}) is 
{\it well-posed} what means that (\ref{deq2}) admits a resolvent $S(t)$.

The so-called resolvent approach to the Volterra equation (\ref{deq2}) has been
introduced many years ago, probably by Friedman and Shinbrot \cite{FrSh},  but
recently has been presented in details in the great monograph by Pr\"uss 
\cite{Pr}. The resolvent approach is a generalization of the semigroup approach.
In the consequence, problems concerning 
convolutions with resolvents (defined below by
(\ref{deq5})) are more difficult than in previous case 
because of lack of the semigroup property.

The main aim of the paper is to provide some introductory 
results for stochastic
convolutions with resolvent operators, analogous to that obtained in e.g.\
\cite{DaPr}, \cite{Tu}
and \cite{DaPrZa}, that is, to extend semigroup approach for our, 
non-semigroup case. 

In order to make the paper self-contained, we formulate definitions and 
auxiliary lemmas necessary for understanding the main results.

Assume that $(\Omega,\mathcal{F},P)$ is a probability space equipped with an 
increasing family of $\sigma$-fields $(\mathcal{F}_t),~t\in [0,T]$
called {\em filtration}. We shall denote by
$\mathcal{F}_{t^+}$ the intersection of all $\sigma$-fields $\mathcal{F}_s$,
$s>t$. We say that filtration is {\em normal} if $\mathcal{F}_0$ contains 
all sets $B\in\mathcal{F}$ with measure $P(B)=0$ and if 
$\mathcal{F}_t=\mathcal{F}_{t^+}$ for any $t\in I$, that is, the filtration is 
right continuous. 

In the paper we assume that filtration 
$(\mathcal{F}_t)_{t\in I}$ is normal. This assumption enables to choose 
modifications of considered stochastic processes with required measurable 
properties.

We will use the following well-known result, see e.g.\ 
\cite{DaPrZa}.

\begin{prop} \label{dprop1}
Let $X(t),~t\in [0,T]$, be a stochastically  
continuous and adapted process with values
in $H$. Then $X$ has a progressively measurable modification.
\end{prop}

In the paper stochastic processes $\Psi$ and $W$ are defined 
 as follows. We consider two separable Hilbert spaces 
$H$ and $U$ and a Wiener
process $W$ on $(\Omega,\mathcal{F},(\mathcal{F}_t)_{t\geq 0},P)$, 
having values in some superspace of $U$ with the nonnegative 
covariance operator $Q\in L(U)$. 
(By $L(U,H)$, $L(U)$ we denote spaces of linear bounded operators 
from $U$ into $H$ and in $U$, respectively.)
We assume that the process $W$ is a
cylindrical one, that is, we do not assume that $\mathrm{Tr}\,Q<+\infty$.
In this case, the process $W$ has continuous paths in some other 
Hilbert space (for
details, see \cite{DaPrZa} or \cite{Ka1}).
Assume that there exists a complete orthonormal set $\{e_k\}\subset U$
of eigenvectors of the operator $Q$ with corresponding 
eigenfunctions $\lambda_k$,
$k=1,2,\ldots;$ so $\mathrm{Tr}Q=\sum_{k=1}^\infty \lambda_k$.
We shall use the following expansion of the process 
$W(t)=\sum_{k=1}^\infty e_k \beta_k(t)$, where $\beta_k(t)$
are independent real Wiener processes with 
$\mathbb{E}(\beta_k^2(t))=\lambda_k t$.
We will need the subspace $U_0:=Q^{1/2}(U)$ of the space $U$, which 
endowed with the inner product 
$\langle u,v\rangle_{U_0}:=\langle Q^{-1/2}u, Q^{-1/2}v\rangle_U$
forms a Hilbert space. Here and in the whole paper we write explicitely
indexes indicating the appropriate space in norms $|\cdot |_{(\cdot)}$
and inner products $\langle\cdot ,\cdot\rangle_{(\cdot)}$ .

This is apparently well-known fact that the construction of the stochastic 
integral with respect to cylindrical Wiener process requires some 
particular terms. Among others, an important role is played by the space of 
Hilbert-Schmidt operators. A linear, bounded operator $C$ acting from $U_0$
into $H$ is called a {\em Hilbert-Schmidt} if
$\sum_{k=1}^{+\infty} |Cu_k|_H^2 < +\infty $, where $\{u_k\}\subset U_0$
is an orthonormal base in $U_0$. The set $L_2(U_0,H)$ of all Hilbert-Schmidt 
operators from $U_0$ into $H$, equipped with the norm
$|C|_{L_2(U_0,H)}:=(\sum_{k=1}^{+\infty} |Cu_k|_H^2)^{1/2}$, is  a separable
Hilbert space. For abbreviation we denote  $L_2^0:=L_2(U_0,H)$.
(For more details concerning that space we refer to \cite{Ba} or \cite{DaPrZa}.)

Let $\Phi(t),~t\in[0,T]$, be a measurable $L_2^0$-valued process. 
We introduce the norms
$$
 ||\Phi||_t := \left\{\mathbb{E}\left( \int_0^t |\Phi(\tau)|_{L_2^0}^2\,d\tau 
 \right) \right\}^{\frac{1}{2}} = \left\{\mathbb{E} \int_0^t 
 \left[ \mathrm{Tr} (\Phi(\tau)Q^{\frac{1}{2}}) (\Phi(\tau)Q^{\frac{1}{2}})^*
 \right] d\tau \right\}^{\frac{1}{2}}, ~~~t\in [0,T].
$$

 By $\mathcal{N}^2(0,T;L_2^0)$ we shall denote a Hilbert space of 
all $L_2^0$-predictable processes $\Phi$ such that $||\Phi ||_T <+\infty$. 

According to the 
theory of stochastic integral with respect to cylindrical Wiener
process (see \cite{DaPrZa} or \cite{Ka1})
we have to assume that $\Psi$ belongs to the space
$\mathcal{N}^2(0,T;L_2^0)$. There is possible to consider a more general
class of integrands, that is, the 
class of $L_2^0$-predictable processes satisfying condition
$
 P\left(\int_0^T | \Psi(\tau) |_{L_2^0}^2 \,d\tau < +\infty\right) = 1.
$
Such processes are called {\em stochastically integrable} on $[0,T]$
and create a linear space denoted by $\mathcal{N}(0,T;L_2^0)$. 
But, in our opinion, it is not worthwhile to study the general case, 
because this assumption makes all formulations of results much more complicated.
Moreover, it produces a new level of difficulty additionally to problems
related to long time memory of the system.

In the whole paper we shall use the following {\sc Volterra Assumptions}
(abbr. ({\sc VA})):
\begin{enumerate}
\item $A:D(A)\subset H\rightarrow H$, is a closed linear unbounded operator
with the dense domain;
\item $a\in L_\mathrm{loc}^1(\mathbb{R}_+)$;
\item the equation (\ref{deq2}) is well-posed and
$S(t),~t\geq 0$, are resolvent operators for the Volterra equation
(\ref{deq2}) determined by the operator $A$ and the function $a$.
\end{enumerate}

For $h\in D(A)$ we define the graph norm as follows:
$|h|_{D(A)}:= (|h|_H^2+|Ah|_H^2)^{\frac{1}{2}}$.
Because $H$ is a separable Hilbert space and $A$ is a closed operator, the
space $(D(A),|\cdot|_{D(A)})$ is a separable Hilbert space.
 
Moreover, we shall study the equation (\ref{deq1}) under the following
{\sc Probability Assumptions} (abbr. ({\sc PA})):
\begin{enumerate}
\item $X_0$ is an $H$-valued, $\mathcal{F}_0$-measurable random variable;
\item $\Psi$ belongs to the space $\mathcal{N}^2(0,T;L_2^0)$, where the finite
interval $[0,T]$ is fixed.
\end{enumerate}

Now, we introduce the definitions of solutions to the stochastic Volterra 
equation (\ref{deq1}). 

\begin{defin} \label{def5}
Assume that conditions (VA) and (PA) hold. 
An $H$-valued predictable process $X(t),~t\in [0,T]$, is said to be a 
{\tt strong solution} to  (\ref{deq1}), if $X$ 
has a version such that $P(X(t)\in D(A))=1$ for almost all $t\in [0,T]$;
for any $t\in [0,T]$, $\int_0^t |a(t-\tau)AX(\tau)|_H d\tau<+\infty$, $P$-a.s.\ and for any 
$t\in [0,T]$ the equation (\ref{deq1}) holds $P$-a.s.
\end{defin}

\noindent{\bf Comment:} 
Because the integral $\int_0^\bullet \Psi(\tau)\,dW(\tau)$ is a continuous 
$H$-valued process
then the above definition yields continuity of the strong solution.

Let $A^*$ denote the adjoint of the operator $A$, with dense domain 
$D(A^*)\subset H$ and the graph norm $|\cdot |_{D(A^*)}$ defined as follows:
$|h|_{D(A^*)}:= (|h|_H^2+|A^*h|_H^2)^{\frac{1}{2}}$, for $h\in D(A^*)$. The
space $(D(A^*),|\cdot|_{D(A^*)})$ is a separable Hilbert space.

\begin{defin} \label{def6}
Let conditions (VA) and (PA) hold.
An $H$-valued predictable process $X(t),~t\in [0,T]$, is said to be a 
{\tt weak solution} to  (\ref{deq1}), if 
$P(\int_0^t|a(t-\tau)X(\tau)|_H d\tau<+\infty)=1$
and if for all $\xi\in D(A^*)$ and all 
$t\in [0,T]$ the following equation holds
$$
 \langle X(t),\xi\rangle_H = \langle X_0,\xi\rangle_H +
 \langle \int_0^t a(t-\tau)X(\tau)\,d\tau,
 A^*\xi\rangle_H +  \langle \int_0^t \Psi(\tau)\,dW(\tau),\xi\rangle_H, 
 \quad P\mathrm{-a.s.}
$$
\end{defin}

\begin{defin} \label{def7}
Assume that (VA) are satisfied and $X_0$ is an $H$-valued 
$\mathcal{F}_0$-measurable random variable. 
An $H$-valued predictable process $X(t),~t\in [0,T]$, is said to be a 
{\tt mild solution} to the stochastic Volterra equation (\ref{deq1}), if
\begin{equation}\label{deq40}
\mathbb{E}\left(\int_0^t |S(t-\tau)\Psi(\tau)|_{L_2^0}^2 d\tau \right)
< +\infty \quad for \quad t\leq T
\end{equation}
and, for arbitrary $t\in [0,T]$,
\begin{equation}\label{deq4}
 X(t) = S(t)X_0 + \int_0^t S(t-\tau)\Psi(\tau)\,dW(\tau),\quad P-a.s.
\end{equation}
\end{defin}

In some cases weak solutions to the equation (\ref{deq1}) 
coincide with mild solutions 
to (\ref{deq1}). In consequence, having results for the convolution 
\begin{equation}\label{deq5}
 W^\Psi(t) :=  \int_0^t S(t-\tau)\Psi(\tau)\,dW(\tau),\quad t\in [0,T],
\end{equation}
where $S(t)$ and $\Psi$ are the same as in (\ref{deq4}),
we obtain results for weak solution to (\ref{deq1}).

\section{Introductory results}

In this section we collect some basic properties of the stochastic convolution
 of the form
\begin{equation} \label{deq6} 
 W^B(t) := \int_0^t S(t-\tau)B\,dW(\tau)
\end{equation} 
in the case when $B\in L(U,H)$.

\begin{lem} \label{dl1} 
 Assume that the operators $S(t),~t\geq 0$, and $B$ are as above,
 $S^*(t), B^*$ are their adjoints, and
\begin{equation} \label{deq7} 
  \int_0^T |S(\tau)B|_{L_2^0}^2\,d\tau =  \int_0^T \mathrm{Tr} 
  [S(\tau) BQB^*S^*(\tau)]\,d\tau < +\infty.
\end{equation} 
Then we have:
\begin{description}
 \item[~~(i)] the process $W^B$ is Gaussian, mean-square continuous on [0,T] 
 and then has a predictable version;
 \item[~(ii)] 
\begin{equation} \label{deq8} 
  \mathrm{Cov}~W^B(t) =  \int_0^t [S(\tau) BQB^*S^*(\tau)]\,d\tau, \quad 
  t\in [0,T];
\end{equation} 
 \item[(iii)] trajectories of the process $W^B$ are P-a.s.\ square integrable
 on [0,T].
\end{description}
\end{lem}
\proof{\begin{description}
 \item[~~(i)] Gaussianity of the process $W^B$ follows from the definition and
 properties of stochastic integral. 
 Let us fix $0\leq t < t+h \leq T$. Then
$$ W^B(t+h)-W^B(t) = \int_0^t [S(t+h-\tau)-S(t-\tau)]BdW(\tau) +
 \int_t^{t+h} \!\!  S(t+h-\tau)BdW(\tau).
$$
Let us note that the above integrals are stochastically independent.
Using the extension of the process $W$ (mentioned in section \ref{defnot})
and properties of stochastic integral with respect to real Wiener processes
(see, e.g., \cite{Ichi}), we have
\begin{eqnarray*} 
 \mathbb{E} | W^B(t+h)-W^B(t) |_H^2 & =& \sum_{k=1}^{+\infty} \lambda_k
 \int_0^t |[ S(t+h-\tau)-S(t-\tau)]Be_k |_H^2 \,d\tau \\
 & + & \sum_{k=1}^{+\infty} \lambda_k \int_t^{t+h} 
 | S(t+h-\tau)Be_k|_H^2 \,d\tau\\
 & := & I_1(t,h) + I_2(t,h) \;.
\end{eqnarray*}
Then, invoking (\ref{deq3a}), the strong continuity of $S(t)$ and the 
Lebesgue dominated convergence theorem, we can pass in $I_1(t,h)$ with 
$h\rightarrow 0$ under the sum and integral signs. 
Hence, we obtain $I_1(t,h) \to 0$ as $h\to 0$.  

Observe that
$$
 I_2(t,h) = \int_t^{t+h} || S(t+h-\tau)BQ^{\frac{1}{2}} ||_{HS}^2 \,d\tau \;, 
$$
where $||\cdot||_{HS}$ denotes the norm of Hilbert-Schmidt operator.
By the condition (\ref{deq7}) we have 
$$
 \int_0^T || S(t) B Q^{\frac{1}{2}} ||_{HS}^2 \,dt < +\infty \;,
$$
what follows that $\lim_{h\rightarrow 0} I_2(t,h) =0$.

The proof for the case $0 \leq t-h < t \leq T$ is similar.
Existence of a predictable version is a consequence of the above continuity and 
Proposition 2. 
\item[~(ii)] Covariance (\ref{deq8}) follows from theory of stochastic integral.
\item[(iii)] From the definition (\ref{deq6}) and assumption (\ref{deq7}) we have 
 the following estimate
\begin{eqnarray*}
\mathbb{E} \int_0^T |W^B(\tau)|_H^2\,d\tau =
  \int_0^T \mathbb{E}|W^B(\tau)|_H^2\,d\tau & = &\\
 = \int_0^T \mathbb{E}\, \left|\int_0^\tau S(\tau-r)BdW(r)\right|_H^2\,d\tau  
 & = & \int_0^T\int_0^\tau 
 |S(r) B|_{L_2^0}^2 \,dr\,d\tau < +\infty\,.
\end{eqnarray*}
Hence, the function $W^B(\cdot)$ may be regarded like random variable with values
in the space $L^2(0,T;H)$.
\hfill $\blacksquare$
\end{description}
}

\noindent{\bf Comment:} Let us emphasize that Cl\'ement and Da Prato
(see \cite{ClDaPr1} and \cite{ClDaPr2}) 
obtained H\"older\-iani\-ty of the trajectories of the 
stochastic convolutions $W_{A,a}(t):=\int_0^t S(t-\tau)dW(\tau)$ in the
case when $A$ is a self-adjoint negative operator in $H$ fulfilling 
some technical
assumptions and when $a\in L_{\mathrm{loc}}^1(\mathbb{R}_+)$ is a 
completely positive function. In that case the operator norm 
$||S(t)||\in [0,1]$
for any $t\in [0,T]$.\\

Analogously like in the theory of evolution equation we can obtain 
the following result.

\begin{theor} \label{t1}
 Assume that the operators $S(t),~t\geq 0$, and $B$ are as above and the 
 condition (\ref{deq7}) holds. Let $X_0$ be a $\mathcal{F}_0$-measurable
 random variable with values in $D(A)$.
 Then the stochastic  Volterra equation (\ref{deq1})  has exactly one weak
 solution which is given by the formula 
 $$ X(t) = S(t) X_0 +\int_0^t S(t-\tau)BdW(\tau). $$
\end{theor}

Now, we formulate an auxiliary result which will be used in the next section.

\begin{lem} \label{dl2}
~Let\/ {\sc Volterra assumptions} hold with the function 
$a \in \,W^{1,1}(\mathbb{R}_+)$. 
Assume that $X$ is a weak solution to (\ref{deq1}) in the case when 
 $\Psi(t)=B$, where $B\in L(U,H)$ and trajectories of $X$
are integrable w.p.\ 1 on $[0,T]$. Then, for any function 
$\xi\in C^1([0,t];D(A^*))$, $t\in [0,T]$, the following formula holds
\begin{eqnarray} \label{deq10} 
  \langle X(t),\xi(t)\rangle_H & = & \langle X_0, \xi(0) \rangle_H +
  \int_0^t  \langle (\dot{a}\star X)(\tau)
  + a(0)X(\tau),A^*\xi(\tau)\rangle_H  d\tau  \nonumber \\ & + &
  \int_0^t \langle \xi(\tau),BdW(\tau)\rangle_H  
  +\int_0^t \langle X(\tau),\dot{\xi}(\tau)\rangle_H d\tau,
\end{eqnarray} 
where dots above $a$ and $\xi$ mean time derivatives and $\;\star\,$ means the
convolution.
\end{lem}

\proof{First, we consider functions of the form 
$\xi(\tau):=\xi_0\varphi(\tau)$,
$\tau\in [0,T]$, where $\xi_0\in D(A^*)$ and $\varphi\in C^1[0,T]$.
For simplicity we omit index $_H$ in the inner product. 
Let us denote 
$ F_{\xi_0}(t) := \langle X(t),\xi_0 \rangle,~ t\in[0,T]. $
 
Using It\^o's formula to the process $ F_{\xi_0}(t)\varphi(t)$, we have
\begin{equation} \label{deq11} 
 d[F_{\xi_0}(t)\varphi(t)] = \varphi(t)dF_{\xi_0}(t) 
 + \dot{\varphi}(t) F_{\xi_0}(t)dt, \quad\quad t\in[0,T].
\end{equation} 
Because $X$ is weak solution to (\ref{deq1}), we have 
\begin{eqnarray} \label{deq12} 
 dF_{\xi_0}(t) & = &\langle\int_0^t \dot{a}(t-\tau)X(\tau)d\tau + a(0)X(t),
 A^*\xi_0\rangle dt +\langle BdW(t),\xi_0\rangle \nonumber \\
 & = & \langle (\dot{a}\star X)(t) +a(0)X(t),A^*\xi_0\rangle dt + 
 \langle BdW(t),\xi_0 \rangle .
\end{eqnarray} 
From (\ref{deq11}) and (\ref{deq12}), we obtain
\begin{eqnarray*}  
 F_{\xi_0}(t)\varphi(t) & = &  F_{\xi_0}(0)\varphi(0) +
 \int_0^t \varphi(s) \langle (\dot{a}\star X)(s) 
 + a(0)X(s), A^*\xi_0\rangle ds \\ &&  + \int_0^t 
 \langle\varphi(s)BdW(s),\xi_0\rangle + \int_0^t \dot{\varphi}(s)
 \langle X(s), \xi_0\rangle ds \\
  & = & \langle X_0, \xi(0) \rangle_H + 
  \int_0^t \langle (\dot{a}\star X)(s)
 + a(0)X(s), A^*\xi(s)\rangle ds \\ 
 && + \int_0^t \langle BdW(s),\xi(s)\rangle 
 +\int_0^t \langle X(s),\dot{\xi}(s)\rangle ds .
\end{eqnarray*} 
Hence, we proved the formula (\ref{deq10}) for functions $\xi$ of the form 
$\xi(s)=\xi_0\varphi(s)$, $s\in [0,T]$. Because such functions form a dense 
subspace in the space  $C^1([0,T];D(A^*))$, the lemma is true.
\hfill $\blacksquare$}\\

\section{Properties in general case}

In this section we consider weak and mild solutions to the equation 
(\ref{deq1}). 

First we study the stochastic convolution defined by (\ref{deq5}),
that is, 
$$
 W^\Psi(t) :=  \int_0^t S(t-\tau)\Psi(\tau)\,dW(\tau),\quad t\in [0,T].
$$

\begin{prop} \label{pr3}
Assume that $S(t),~t\geq 0$, are (as earlier) the resolvent operators corresponding
to the Volterra equation (\ref{deq2}). Then, for arbitrary process 
$\Psi\in\mathcal{N}^2(0,T;L_2^0)$, the process $W^\Psi(t),~t\geq 0$, given by 
(\ref{deq5}) has a predictable version.
\end{prop}

\proof{Because proof of Proposition \ref{pr3} is analogous to some schemes 
in theory of stochastic integral (see,e.g., \cite[Chapter 4]{LiSh}) 
we provide only an outline of proof.

First, let us notice that the process $S(t-\tau)\Psi(\tau)$, where 
$\tau\in[0,t]$,
belongs to $\mathcal{N}^2(0,T;L_2^0)$, because
$\Psi\in\mathcal{N}^2(0,T;L_2^0)$.
Then we may use the apparently well-known estimate
(see, e.g., Proposition 4.16 in \cite{DaPrZa}):  
for arbitrary $a>0,~b>0$ and $t\in[0,T]$ 
\begin{equation} \label{g1}
P(|W^\Psi(t)|_H>a) \leq \frac{b}{a^2}+P\left( \int_0^t  
 |S(t-\tau)\Psi(\tau)|_{L_2^0}^2 d\tau >b\right)\;.
\end{equation}
Because the resolvent operators $S(t),~t\geq 0$, are uniformly bounded on compact
itervals (see \cite{Pr}), there exists a constant $C>0$ such that $||S(t)||\leq C$
for $t\in[0,T]$. So, we have 
$|S(t-\tau)\Psi(\tau)|_{L_2^0}^2\leq C^2|\Psi(\tau)|_{L_2^0}^2$, $\tau\in[0,T]$.

Then the estimate (\ref{g1}) may be rewritten as 
\begin{equation} \label{g2}
P(|W^\Psi(t)|_H>a) \leq \frac{b}{a^2}+P\left( \int_0^t  
|\Psi(\tau)|_{L_2^0}^2 d\tau >\frac{b}{C^2}\right)\;.
\end{equation}

Let us consider predictability of the process $W^\Psi$ in two steps. 
In the first step we assume that $\Psi$ is an elementary process
understood in the sense given in section 4.2 in \cite{DaPrZa}. In this
case the process $W^\Psi$ has a predictable version by Lemma \ref{dl1}, 
part (i).

In the second step $\Psi$ is an arbitrary process belonging to 
$\mathcal{N}^2(0,T;L_2^0)$. Since elementary processes form a dense 
set in the space $\mathcal{N}^2(0,T;L_2^0)$, there exists a sequence 
$(\Psi_n)$ of elementary processes such that for arbitrary $c>0$
\begin{equation} \label{g3}
P\left( \int_0^T |\Psi(\tau)-\Psi_n(\tau)|_{L_2^0}^2 d\tau 
 >c\right)\stackrel{n\rightarrow +\infty}{\longrightarrow} 0\;.
\end{equation}

By the previous part of the proof the sequence $W_n^\Psi$ of convolutions
$$ W_n^\Psi (t):= \int_0^t S(t-\tau)\Psi_n(\tau)dW(\tau)$$
converges in probability. Hence, it has a subsequence converging almost
surely. This implies the predictability of the convolution 
$W^\Psi (t), ~t\in [0,T]$.
~\hfill $\blacksquare$}\\

\begin{prop} \label{pr3a}
Assume that $\Psi\in\mathcal{N}^2(0,T;L_2^0)$. Then the process 
$W^\Psi(t),~t\geq 0$, defined by (\ref{deq5}) has square integrable 
trajectories.
\end{prop}

\proof{We have to prove that $\mathbb{E}\int_0^T |W^\Psi(t)|_H^2 dt<+\infty$.
From Fubini's theorem and properties of stochastic integral
\begin{eqnarray*}
 \mathbb{E}\int_0^T \left| \int_0^t S(t-\tau)\Psi(\tau)dW(\tau)\right|_H^2 dt
  & =& \int_0^T \left[\mathbb{E}\left| 
 \int_0^t S(t-\tau)\Psi(\tau)dW(\tau)\right|_H^2\right] dt \\
 = \int_0^T\int_0^t |S(t-\tau)\Psi(\tau)|_{L_2^0}^2\; d\tau dt & \leq &
 M \int_0^T\int_0^t |\Psi(\tau)|_{L_2^0}^2 \; d\tau dt  < +\infty.  \\
  \mbox{~(from~boundness~of operators~} S(t) 
  &&
  \mbox{and~because~} \Psi(\tau) \mbox{~are~Hilbert-Schmidt)} 
\end{eqnarray*}
\hfill $\blacksquare$}

In the below result, the notions "parabolic" and "3-monotone" are understood 
in the sense defined by Pr\"uss \cite[Section 3]{Pr}.

\begin{prop} \label{pr4}
Assume that (\ref{deq1}) is parabolic, (VA) are satisfied and 
the kernel function  $a$ is 3-monotone. 
Let $X$ be a predictable process with integrable trajectories.
Assume that $X$ has a version such that $P(X(t)\in D(A))=1$ for almost all
$t\in [0,T]$ and (\ref{deq40}) holds.
If for any $t\in [0,T]$ and $\xi\in D(A^*)$
\begin{equation}\label{deq9}
 \langle X(t),\xi\rangle_H = \langle X_0,\xi\rangle_H +
 \int_0^t \langle 
 a(t-\tau)X(\tau),A^*\xi\rangle_H d\tau 
 + \int_0^t \langle\xi,\Psi(\tau) dW(\tau)\rangle_H,  ~~P-a.s.,
\end{equation}
then
\begin{equation}\label{deq9a}
X(t) = S(t)X_0 +
 \int_0^t S(t-\tau) \Psi (\tau) dW(\tau), \quad t\in[0,T].
\end{equation}
\end{prop}

\proof{For simplicity we omit index $_H$ in the inner product. 
First, we see, analogously like in Lemma \ref{dl2}, 
that if (\ref{deq9}) is satisfied, then 
\begin{eqnarray}\label{deq14}
 \langle X(t),\xi(t)\rangle &=& \langle X_0,\xi(0)\rangle +
 \int_0^t \langle (\dot{a}\star X)(\tau) +
 a(0)X(\tau),A^*\xi(\tau)\rangle d\tau \nonumber \\
 &+& \int_0^t \langle \Psi(\tau) dW(\tau),\xi(\tau)\rangle  
 + \int_0^t \langle X(\tau),\dot{\xi}(\tau)  \rangle d\tau, \quad
 \mathrm{~P-a.s.} 
\end{eqnarray}
holds for any $\xi\in C^1([0,t],D(A^*))$ and $t\in [0,T]$.

Because (\ref{deq1}) is parabolic and $a$ is 2-regular (what is implied 
by 3-monotone), then, by \cite[Theorem 3.1]{Pr},
there exists a resolvent $S\in C^1((0,+\infty);L(H))$ for (\ref{deq1}).

Now, let us take $\xi(\tau):=S^*(t-\tau)\zeta$ with $\zeta\in D(A^*)$, 
 $\tau\in [0,t]$. 
The equation (\ref{deq14}) may be written like
\begin{eqnarray*}  
 \langle X(t),S^*(0)\zeta\rangle & = & \langle X_0,S^*(t)\zeta\rangle
 + \int_0^t \langle (\dot{a}\star X)(\tau)
 + a(0)X(\tau),A^*S^*(t-\tau)\zeta\rangle d\tau \nonumber \\
  &+& \int_0^t \langle \Psi(\tau) dW(\tau),S^*(t-\tau)\zeta\rangle
 + \int_0^t \langle X(\tau),(S^*(t-\tau)\zeta)'\rangle d\tau, 
\end{eqnarray*}
where derivative ()' in the last term is taken over $\tau$. 

Next, using $S^*(0)=I$, we rewrite
\begin{eqnarray} \label{deq15}
 \langle X(t),\zeta\rangle &=& \langle S(t)X_0,\zeta\rangle +
 \int_0^t \langle S(t-\tau)A\left[\int_0^\tau
 \dot{a}(\tau-\sigma)X(\sigma)d\sigma +a(0)X(\tau)\right],\zeta\rangle 
 d\tau \nonumber\\
   &+& \int_0^t \langle S(t-\tau)\Psi(\tau)dW(\tau),\zeta\rangle +
   \int_0^t \langle \dot{S}(t-\tau)X(\tau),\zeta\rangle d\tau. 
\end{eqnarray}
To prove  (\ref{deq9a}) it is enough to show that the sum of the first 
integral and the third one in the equation (\ref{deq15}) gives zero.

We use properties of resolvent operators and the derivative $\dot{S}(t-\tau)$
with respect to~$\tau$. Then
\begin{eqnarray*}  
 I &:=& \left\langle \int_0^t \dot{S}(t-\tau)X(\tau) d\tau,\zeta \right\rangle =
  \left\langle -\int_0^t \dot{S}(\tau)X(t-\tau) d\tau ,\zeta \right\rangle \\
 &=& \left\langle -\left(
  \int_0^t\left[ \int_0^\tau \dot{a}(\tau-s)AS(s)ds\right] X(t-\tau) d\tau 
 - \int_0^t a(0)AS(\tau) X(t-\tau) d\tau \right) ,\zeta \right\rangle \\[2mm]
 & = & \langle -([A(\dot{a}\star S)(\tau) \star X](t)+a(0)A(S\star X)(t)) 
 ,\zeta\rangle .
\end{eqnarray*}
The kernel function $a$ is 3-monotone, so $a\in C^1(\mathbb{R}_+)$, and then
has bounded variation. Hence, the convolution 
$(a\star S)(\tau)$ has sense (see \cite[Section 1.6]{Pr} or \cite{ArKe}).

Since
$$ \int_0^t \langle a(0)AS(t-\tau)X(\tau),\zeta\rangle d\tau =
  \int_0^t \langle a(0)AS(\tau)X(t-\tau),\zeta\rangle d\tau $$
 and
\begin{eqnarray*} 
J &:=& \int_0^t \langle S(t-\tau)A\left[\int_0^\tau
 \dot{a}(\tau-\sigma)X(\sigma)d\sigma\right],\zeta\rangle d\tau =
 \int_0^t \langle AS(t-\tau)(\dot{a}\star X)(\tau),\zeta\rangle d\tau = \\[2.5mm]
 &\!\!=\!\!& \langle A(S\star (\dot{a}\star X)(\tau))(t),\zeta\rangle =
 \langle A((S\star\dot{a})(\tau)\star X)(t),\zeta\rangle
 \quad \mbox{for~any~~} \zeta\in D(A^*),
\end{eqnarray*} 
so $J=-I$, hence $J+I=0$. 
This means that (\ref{deq9a}) holds for any $\zeta\in D(A^*)$. Since 
$D(A^*)$ is dense in $H^*$, then (\ref{deq9a}) holds.
\hfill $\blacksquare$}

\noindent{\bf Remark:}
In Proposition \ref{pr4}, the assumption that $a$ is 3-monotone may be replaced
by both: 2-regularity of $a$ and $a\in BV_\mathrm{loc}(\mathbb{R}_+)$.  

\noindent{\bf Comment:} 
Proposition \ref{pr4} shows that under particular conditions a weak 
solution to (\ref{deq1}) is a mild solution to the equation (\ref{deq1}).

\begin{prop} \label{pr5}
Let {\sc Volterra assumptions} be satisfied. 
If $\Psi\in\mathcal{N}^2(0,T;L_2^0)$ and $\Psi(\cdot,\cdot)(U_0)\subset D(A),$
 P--a.s., then the stochastic convolution 
$W^\Psi$ fulfills the equation (\ref{deq9}) with $X_0\equiv 0$.
\end{prop}
\proof{Let us notice that the process $W^\Psi$ has integrable trajectories.
For any $\xi\in D(A^*)$ we have
\begin{eqnarray*}
 \int_0^t \langle a(t-\tau)W^\Psi(\tau),A^*\xi\rangle_Hd\tau &\equiv& 
 ~\mbox{(from~ (\ref{deq5}))} \\
 &\equiv& \int_0^t \langle a(t-\tau) \int_0^\tau S(\tau-\sigma)\Psi(\sigma) 
 dW(\sigma),A^*\xi\rangle_Hd\tau =\\
  (\mbox{from~Dirichlet's~formula}& \mbox{and} & 
  \mbox{stochastic Fubini's theorem})  \\
 &=& \int_0^t \langle\left[\int_\sigma^t a(t-\tau)S(\tau-\sigma)d\tau\right]
 \Psi(\sigma)  dW(\sigma),A^*\xi\rangle_H\\
 &=& \langle \!\int_0^t\! \left[\int_0^{t-\sigma}\! a(t-\sigma-z)S(z)dz\right]
 \!\Psi(\sigma) dW(\sigma),A^*\xi\rangle_H\\
 (\mbox{where~} z:=\tau-\sigma&& \mbox{~and~from definition~of~convolution})\\
  &=& \langle \int_0^t A [ (a\star S)(t-\sigma)]\Psi(\sigma) dW(\sigma),
  \xi\rangle_H = \\
  (\mbox{from~the~resolvent~equation~(3),} &&  \mbox{because~~}
     A(a\star S)(t-\sigma)x = (S(t-\sigma)-I)x, \\
    && ~\mbox{~where~} x\in D(A))  \\
 &=& \langle \int_0^t [S(t-\sigma)-I]\Psi(\sigma) dW(\sigma),\xi\rangle_H = \\
 = \langle \int_0^t S(t-\sigma)\Psi(\sigma) dW(\sigma),\xi\rangle_H
 &-&\langle \int_0^t \Psi(\sigma) dW(\sigma),\xi\rangle_H .
\end{eqnarray*}
Hence, we obtained the following equation
$$
 \langle W^\Psi(t),\xi\rangle_H = \int_0^t \langle a(t-\tau)W^\Psi(\tau), 
 A^*\xi\rangle_Hd\tau + \int_0^t \langle \xi,\Psi(\tau)dW(\tau)\rangle_H
$$
for any $\xi\in D(A^*)$.
\hfill $\blacksquare$}

\begin{cor} \label{c1}
Let {\sc Volterra assumptions} hold with a bounded operator $A$. 
If $\;\Psi$ belongs to $\mathcal{N}^2(0,T;L_2^0)$  then
\begin{equation} \label{deq16}
 W^\Psi(t) = \int_0^t a(t-\tau)AW^\Psi(\tau)d\tau 
 + \int_0^t \Psi(\tau) dW(\tau)\,.
\end{equation}
\end{cor}

\noindent{\bf Comment:} 
The formula (\ref{deq16}) says that the convolution $W^\Psi$ is a strong
solution to (\ref{deq1}) with $X_0\equiv 0$ if the operator $A$ is bounded.

The below theorem is a consequence of the results obtained up to now.

\begin{theor} \label{t3}
Suppose that (VA) and (PA) hold. 
Then a strong solution (if exists) is always a weak solution of (\ref{deq1}). 
If, additionally, assumptions of Proposition \ref{pr4} are satisfied, 
a weak solution  is a mild solution to the 
Volterra equation (\ref{deq1}). Conversely, under conditions of Proposition 
\ref{pr5}, if $X_0\equiv 0$ a mild solution $X$ is also a weak solution to 
(\ref{deq1}).
\end{theor}

\section{Some estimates}

In this section we provide two estimates for stochastic convolution 
(\ref{deq5}). 
Some considerations (see, e.g.\ \cite{Ka2}, where maximal type inequalities
for the equation (\ref{deq1}) were studied) show that in general case, that is 
when (VA) are supposed only, it is very difficult 
to say something interesting about regularity of the convolution (\ref{deq5}).
Similar situation was in the semigroup case, see \cite{DaPrZa}, 
where regularity results have been received under additional assumptions on
semigroups, for instance when contractions or analytical semigroups were studied.

As we have already written, Cl\'ement and Da Prato (see, \cite{ClDaPr1} and
\cite{ClDaPr2}) obtained some regularity results in the case when $A$
was a self-adjoint operator satisfying some technical assumptions, the function
$a$ was completely positive and when $\Psi(t)\equiv I$, that is for the 
convolution $\widetilde{W}(t):=\int_0^t\,S(t-\tau)dW(\tau)$. In their case,
$||S(t)||\leq 1$, what is a resolvent analogon of contraction semigroup.

\begin{theor} \label{t4}
If $\; \Psi\in \mathcal{N}^2(0,T;L_2^0)$ then the following estimate holds
 \begin{equation} \label{deq18}
  \sup_{t\leq T} \;\mathbb{E} (|W^\Psi (t)|_H) \leq  C\,
  M_T\,  \mathbb{E} \left( \int_0^T|\Psi(t)|_{L_2^0}^2\, dt 
  \right)^{\frac{1}{2}},
 \end{equation}
 where $C$ is a constant and $M_T=\sup_{t\leq T} ||S(t)|| $.
\end{theor}

\noindent{\bf Comment:} The estimate (\ref{deq18}) seems to be rather coarse.
It comes directly from the definition of stochastic integral. 
Since (\ref{deq18}) reducest to the Davis inequality for msrtingales if 
$S(\cdot)=I$, the constant $C$ appeared on the right hand side.
Unfortunately, we can 
not use more refined tools, for instance It\^o's formula (see, e.g.\ \cite{Tu}
for Tubaro's estimate), because the process $W^\Psi$ is not enough regular. 

The next theorem is a consequence of Theorem \ref{t4}.

\begin{theor} \label{t5}
 Assume that $\Psi\in\mathcal{N}^2(0,T;L_2^0)$.  Then 
 $$\sup_{t\leq T}\, \mathbb{E}(|W^\Psi(t)|_H)
 \leq \widetilde{C}(T)\,|\Psi|_{\mathcal{N}^2(0,T;L_2^0)},$$
 where a constant $\widetilde{C}(T)$ depends on $T$.
\end{theor}

\proof{From (\ref{deq5}) and property of stochastic integral we have 
\begin{eqnarray*}
 \mathbb{E}(|W^\Psi(t)|_H) = \mathbb{E}\left(|\int_0^t 
 S(t-\tau)\Psi(\tau)dW(\tau)|_H\right) &\leq & 
 C\,\mathbb{E}\left(\int_0^t | S(t-\tau)\Psi(\tau)|_{L_2^0}^2\,d\tau
 \right)^{\frac{1}{2}}  \\
 \mbox{(from~writing~out~the~Hilbert-Schmidt~norm)}&\leq & 
 C\,\mathbb{E}\left(\int_0^t ||S(t-\tau)||^2 \, |\Psi(\tau)|_{L_2^0}^2
 \,d\tau \right)^{\frac{1}{2}} \\
   \leq C\, M_T\, \mathbb{E}\left(\int_0^t |\Psi(\tau)|_{L_2^0}^2
 \,d\tau \right)^{\frac{1}{2}} & \le & (\mbox{by~the~H\"older~inequality})\\
  \le  C\, M_T\, \left( \mathbb{E}\int_0^t |\Psi(\tau)|_{L_2^0}^2\,d\tau 
 \right)^{\frac{1}{2}}
 & = & \widetilde{C}(T)\, |\Psi|_{\mathcal{N}^2(0,T;L_2^0)}\;,
\end{eqnarray*}
where $M_T$ is as above and $\widetilde{C}(T)= C\,M_T$.
\hfill $\blacksquare$}\\

Summing up the paper, it is worth to emphasize that better and more 
sophisticated regularity results should be obtained for exponentially 
bounded and analytical resolvents.
The situation is similar to that for semigroup case, when the best results
are reached for analytical semigroups.\\


\begin{thebibliography}{AB}

\bibitem{ArKe} Arendt W.\ and Kellerman H., \textit{Integral solutions of
Volterra integrodifferential equations and applications} in G.\ Da Prato, M.\
Iannelli, eds.\ - Volterra Integrodifferential Equations in Banach Spaces and
Applications, pp.\ 21-51, Harlow, Essex, 1989, Longman Sci.\ Tech. 

\bibitem{Ba} Balakrishnan A.V.,  \textit{Applied Functional Analysis},
Springer-Verlag, New York, 1981.

\bibitem{ClDaPr1}  Cl\'{e}ment Ph. and Da Prato G., {\it Some results on
  stochastic convolutions arising in Volterra equations perturbed by noise},
  Rend. Math. Acc. Lincei s.9, \textbf{7} (1996), 147-153.
 
\bibitem{ClDaPr2} Cl\'{e}ment Ph. and Da Prato G., {\it White noise
  perturbation of the heat equation in materials with memory}, 
Dynamic Systems and Applications \textbf{6} (1997), 441-460.
  
\bibitem{ClDaPP} Cl\'{e}ment Ph., Da Prato G.\ and Pr\"uss J., {\it White 
noise perturbation of the equations of linear parabolic viscoelasticity},
Rend.\ Inst.\ Mat.\ Univ.\ Trieste {\bf 29} (1997), 207--220. 

\bibitem{DaPr}Da Prato G., {\it Regularity results of a covolution
stochastic integral and applications to parabolic stochastic equations in
Hilbert space}, Conferenze del Seminario Matematico 
dell'Universit\'a di Bari, No. {\bf 182}, Laterza. 
  
\bibitem{DaPrZa} Da Prato G.\ and Zabczyk J., {\it Stochastic
equations in infinite dimensions}, Cambridge University Press,
Cambridge, 1992.

\bibitem{FrSh} Friedman A. and Shinbrot M., {\it Volterra integral equations 
 in Banach space}, Trans.\ Am.\ Math.\ Soc. {\bf 126} (1967) 131--179.

\bibitem{Ichi}  Ichikawa A., {\it Stability of semilinear stochastic evoltion
equations}, J.\ Math.\ Anal.\ Appl., {\bf 90} (1982) 12--44.


\bibitem{Ka1}  Karczewska A., {\it Stochastic integral with respect to 
cylindrical Wiener process}, Annales Universitatis Mariae 
Curie-Sk\l odowska, Vol.\ LII.\ 2,9 (1998) 79--93. 
 http://xxx.lanl.gov/abs/math.PR/0511512

\bibitem{Ka2}  Karczewska A., {\it Maximal type inequalities for linear 
 stochastic Volterra equations}, Int.\ J.\ Pure Appl.\ Math.\ 
 {\bf 24} (2005) no.\ 1, 111--121.
 http://xxx.lanl.gov/abs/math.PR/0412496

\bibitem{LiSh} Liptser R.S and Shiryayev A.N.,
 \textit{Statistics of random processes.II} Applications of Mathematics, Vol.\
 {\bf 6}, Springer, New York, 1973.
  
\bibitem{Pr} Pr\"uss J.,  \textit{Evolutionary integral equations and
applications}, Birkh\"auser, Ba\-sel, 1993.

\bibitem{Tu}Tubaro L., \textit{An estimate of Burkholder type for stochastic
processes defined by the stochastic integral}, Stoch.\ Anal.\  Appl.\, {\bf 2} 
(1984) 187--192.
 
\end{thebibliography}
\end{document}